\documentclass[matzei,11pt,twoside]{amsart} 

\title{Euler homology}
\author{Julia Weber}
\address{Max-Planck-Institut f\"ur Mathematik \newline \indent Vivatsgasse 7 \newline \indent 53111 Bonn \newline \indent Germany}
\email{jweber@mpim-bonn.mpg.de}

\subjclass[2000]{primary: 55N20, 57R90, 57R99; secondary: 55N91, 57R85}
\keywords{generalized homology theory, Euler characteristic, equivariant}

\usepackage[latin1]{inputenc} 
\usepackage{bbm}
\usepackage{amsmath}          
\usepackage{amssymb}
\usepackage{amsthm}
\usepackage{xspace}
\usepackage[all]{xy}\CompileMatrices
\usepackage{amsfonts}
\usepackage{amscd}
\usepackage{hyperref}

\DeclareMathOperator{\PD}{PD}
\DeclareMathOperator{\IIm}{Im}
\DeclareMathOperator{\pr}{pr}
\DeclareMathOperator{\ind}{ind}
\newcommand{\Ab}{\mathcal{A}b}

\newcommand{\GTop}{G\text{-}\mathcal{T}op}

\newcommand{\Z}{\ensuremath{\mathbb{Z}}\xspace}
\newcommand{\R}{\ensuremath{\mathbb{R}}\xspace}
\newcommand{\Nat}{\ensuremath{\mathbb{N}}\xspace}
\newcommand{\Zz}{\ensuremath{\Z/2}\xspace}
\newcommand{\N}{\ensuremath{\mathcal{N}}\xspace}
\newcommand{\Cl}{\ensuremath{\mathcal{C}}\xspace}
\renewcommand{\H}{\ensuremath{\mathcal{H}}\xspace}
\newcommand{\K}{\ensuremath{\mathcal{K}}\xspace}
\newcommand{\oC}{\ensuremath{\text{\emph{\r{C}}}}\xspace}

\newcommand{\oM}{\ensuremath{\text{\emph{\r{M}}}}\xspace}
\newcommand{\oS}{\ensuremath{\text{\emph{\r{S}}}}\xspace}
\newcommand{\oT}{\ensuremath{\text{\emph{\r{T}}}}\xspace}

\newtheorem{Def}{Definition}[section]
\newtheorem{Thm}[Def]{Theorem}
\newtheorem{Cor}[Def]{Corollary}
\newtheorem{Ex}[Def]{Example}
\newtheorem{Prop}[Def]{Proposition}
\newtheorem{Lem}[Def]{Lemma}
\newtheorem*{Rmk}{Remark}

\begin{document}

\maketitle

\begin{abstract}

We geometrically construct a homology theory that generalizes the
Euler characteristic mod 2 to objects in the unoriented cobordism ring
$\N_*(X)$ of a topological space $X$. This homology theory $Eh_*$ has
coefficients $\Zz$ in every nonnegative dimension. There exists a
natural transformation $\N_*(X)\to Eh_*(X)$ that for $X=pt$ assigns to
each smooth manifold its Euler characteristic mod 2. The homology
theory is constructed using cobordism of stratifolds, which are singular objects defined below. 

An isomorphism $Eh_*(X)\cong  H_*(X;\Zz)\otimes_{\Zz} \Zz[t] $ of graded $\N_*$-mod\-ules is shown for any CW-complex $X$. For discrete groups $G$, we also define an equivariant version of the homology theory $Eh_*$, generalizing the equivariant Euler characteristic. 

\end{abstract}

\section{Introduction}

We construct a homology theory using bordism of Euler stratifolds, singular objects defined below. This homology theory $Eh_*$ has coefficients $\Zz$ in every nonnegative dimension. There is a natural transformation from the unoriented bordism theory $\N_*(X)$ to $Eh_*(X)$ obtained by regarding a smooth manifold as an Euler stratifold. For $X=pt$ this assigns to each smooth manifold its Euler characteristic mod 2. 

For CW-complexes $X$, we obtain $Eh_*(X)\cong H_*(X;\Zz)\otimes_{\Zz} \Zz[t]$ with $\mbox{deg } t=1$. The natural transformation $\N_*(X) \to Eh_*(X)$ maps the element $[M,f]\in \N_n(X)$ to the element $\sum_{0 \leq j \leq n} f_* \PD w_{n-j}(M) \otimes t^{n-j} \in H_*(X;\Zz)\otimes_{\Zz} \Zz[t] \cong Eh_*(X)$, where $\PD w_{n-j}$ is the Poincaré dual of the $({n-j})$-th Stiefel-Whitney class. Thus the homology theory $Eh_*$ allows a geometric detection of these bordism invariants.

The objects used in the geometric construction are stratifolds,
introduced by Kreck~\cite{kreck}. Stratifolds are differentiable
singular objects that generalize smooth manifolds. They are very
useful for defining homology theories geometrically. One can define
singular homology as a bordism theory using stratifolds by placing
certain restrictions on them~\cite{kreck}. It is a natural question to
consider other sets of restrictions, to ask whether these also lead to
bordism theories and, if so, to analyze the homology theories thus
obtained. 

The restrictions on the Euler characteristics of certain subobjects
are motivated by~\cite{akin}. The objects considered there are
polyhedra, so the context is quite different. Akin constructs an
ungraded bordism theory isomorphic to total $\Zz$-homology. In contrast to
Akin, we construct a graded homology theory and a graded isomorphism
$Eh_*(X)\cong  H_*(X;\Zz)\otimes_{\Zz} \Zz[t] $ of graded
$\N_*$-mod\-ules. Our theory recovers Akin's when we set $t=1$. (This
corresponds to forgetting the dimension of a stratifold.) Akin
mentions that the coefficients of a graded homology theory would be
$\Zz$ in every nonnegative dimension, but this is not analyzed any
further. An advantage of our approach is that we obtain a
canonical natural transformation $\N_*(X) \to Eh_*(X)$ of graded
$\N_*$-mod\-ules which cannot be constructed using the ungraded
theory. 

The use of these restrictions on the Euler characteristic in the
context of stratifolds was suggested by Kreck. I want to express
sincere thanks for proposing this question, for introducing me to the
fascinating subject of stratifolds, and for many animated discussions
while this work was developed. I also thank the referee for his
detailed report and helpful comments.

\section{Euler stratifolds}

Stratifolds are a generalization of smooth manifolds. They were
introduced by Kreck~\cite{kreck}. Euler stratifolds are defined using
parametrized stratifolds~\cite[Chapter~2.3, Example~9]{kreck}. For the
convenience of the reader we recall definitions and techniques from
this theory. Details can be found in~\cite{kreck}. 

We always consider smooth manifolds in the sequel, thus whenever manifolds are mentioned they are implicitly assumed to be smooth. Pa\-ra\-me\-trized stratifolds are defined using c-manifolds, manifolds with germs of collars.

\begin{Def} \label{defcmap}
A \emph{collar} of a manifold $M$ with boundary $\partial M$ is a
smooth map $\varphi\colon \partial M\times [0,\varepsilon)\to M$
satisfying $\varphi(x,0)=x$ for all $x\in \partial M$ and that is a
diffeomorphism onto its image. Two such collars $\varphi\colon
\partial M\times [0,\varepsilon)\to M$ and $\varphi'\colon \partial
M\times [0,\varepsilon')\to M$are called \emph{equivalent} if there
exists an $\varepsilon'' \leq \min(\varepsilon,\varepsilon')$ such
that the restrictions of $\varphi$ and $\varphi'$ to $ \partial
M\times [0,\varepsilon'')$ coincide.

A \emph{c-manifold} $M$ is a manifold $M$ with (possibly empty)
boundary equipped with an equivalence class of collars. 

A \emph{c-map} between c-manifolds $M$ and $N$ is a smooth map
$f\colon M\to N$ such that there exist representatives
$\varphi_M\colon \partial M\times [0,\varepsilon)\to M$ and
$\varphi_N\colon \partial N\times [0,\varepsilon)\to N$ of the germs
of collars for which either 
\[
f(\varphi_M(x,t))= f(x) \text{ for all }(x,t)\in \partial M\times [0,\varepsilon)
\]
 or 
\[
f(\varphi_M(x,t))=
\varphi_{N}(f(x),t) \text{ for all } (x,t)\in \partial M\times
[0,\varepsilon).
\]
The latter case is only possible if $f(x)\in
\partial N$.
\end{Def}

\begin{Ex}
The manifold $I:=[0,1]$ has a canonical equivalence class of collars
given by the maps $\varphi:=(id,1-id)\colon  (\{0\},\{1\})\times
[0,\varepsilon) \to [0,1]$. We view $I$ as a c-manifold equipped with
this equivalence class of collars if nothing else is mentioned in the
context. The same holds for $[0,1), (0,1]$ and similar spaces.
\end{Ex}

Any manifold can be equipped with a collar. A manifold $M$ without boundary is canonically a c-manifold, equipped with the equivalence class of the collar $\varphi\colon \emptyset \to M$. One defines products of c-manifolds by smoothing of collars. C-manifolds can be seen as building blocks of parametrized stratifolds.

\begin{Def} \label{defstrat}
An $n$-dimensional \emph{parametrized stratifold} $S$ is a locally
compact Hausdorff space $S$ with countable basis together with proper maps $f_0\colon S_0 \to S, \ldots , f_n\colon S_n\to S$ such that 
\begin{itemize}
\item $S_i$ is an $i$-dimensional c-manifold 
\item $f_i|_{\oS_i}$ is a homeomorphism onto its image and $S=\coprod_{i\leq n} f_i(\oS_i)$
\item $f_i(\partial S_i) \subseteq S^{(i-1)}:=\bigcup_{j\leq i-1} f_j(S_j)$
\item $f_i|_{\partial S_i}\colon \partial S_i \to S^{(i-1)}$ is
      ``smooth'': This means that for all continuous maps $\rho\colon S^{(i-1)} \to \R$ such that $\rho f_j\colon S_j\to \R$ is a c-map for all $j\leq i-1$ the composition $\rho f_i|_{\partial S_i}\colon \partial S_i \to \R$ is smooth.
\end{itemize}
We call $S^{i}$ the \emph{$i$-stratum} of $S$ and $S^{(i-1)}$
the \emph{$(i-1)$-skeleton} of $S$. 
\end{Def}

We also call parametrized stratifolds \emph{p-stratifolds}. One can
visualize a p-stratifold by imagining it being built inductively from
the smooth manifolds $S_i$. One starts with $S_0$, which is just a
collection of points, then takes $S_1$, which is a collection of
$1$-spheres and intervals, these intervals being glued to $S_0$ at the
boundary points via $f_1$, then glues $S_2$ to this structure via
$f_2$, and so forth. The equivalence classes of collars and the maps
$f_i$ are explicitly part of the structure. 

Strata can be disconnected. In every dimension from $0$ to $n$, there
is exactly one stratum. Note that a stratum can be empty. This means
that we cannot read off the dimension of a p-stratifold from the
topological space $S$ alone: The topmost strata may be empty, but are
nevertheless part of the structure. 

In our definition, the dimension $n$ is always finite, so we only have
finitely many strata. We could also define infinite-dimensional
stratifolds, but we will not need these in the sequel.

\begin{Ex}
An $m$-dimensional manifold $M$ without boundary can be seen as an $m$-dimensional p-stratifold by taking $S=M$, $S_0=\ldots=S_{m-1}=\emptyset$, and $S_m=M, f_m=id_M$. If the context does not imply anything else, we mean this p-stratifold whenever we view a manifold as a p-stratifold.
\end{Ex}
\begin{Ex}
An $m$-dimensional manifold $M$ can be seen as an $n$-di\-men\-sional p-stratifold for all $n\geq m$ by filling it up with empty strata $S_{m+1}=\emptyset$, $\ldots$, $S_{n}=\emptyset$. We denote this p-stratifold by $M_n$. A special case of this are the p-stratifolds $pt_n$, for $n\geq 0$.  
\end{Ex}
\begin{Ex}
The open cone $\oC S:=S\times (0,1]/S\times \{1\}$ over a p-stratifold $S$ has an induced p-stratifold structure given by  $f_0\colon (\oC S)_0:=pt\to \overline{S\times \{1\}}\subseteq S\times (0,1]/S\times \{1\}$ and $f_{i+1}=f_i\times id_{(0,1]}\colon (\oC S)_{i+1}:=S_i\times (0,1]\to S\times (0,1]/S\times \{1\}$ for all $0\leq i \leq n$. 
\end{Ex}

An open subset $U\subseteq S$ of a p-stratifold $S$ that is compatible with the collar structure is a p-stratifold, with the induced p-stratifold structure given by intersection with $U$:

\begin{Def} \label{defa}
Let $U$ be an open subset of the $n$-dimensional p-stratifold $S$ such that for all $i$ there is a collar $\varphi_{S_i}\colon \partial S_i\times [0,\varepsilon)\to S_i$ in the given equivalence class of collars such that for all $x\in f_i^{-1}(U)\cap \partial S_i$ we have $f_i(\varphi_{S_i}(\{x\}\times [0,\varepsilon)))\in U$. Then we define an \emph{induced p-stratifold structure} on $U$ by $U_i:=S_i\cap f_i^{-1}(U)$ and $(f_U)_i:=(f_S)_i|_{U_i}$ for all $0\leq i \leq n$. For the collars we take the equivalence classes of the collars $\varphi_{S_i}|_{\partial U_i \times [0,\varepsilon)}$.
\end{Def}  

One can also define a product of p-stratifolds:

\begin{Def} \label{defb}
 The \emph{product} $S\times T$ of two p-stratifolds $S$ and $T$ of dimensions $k$ respectively $l$ is given by the topological space $S\times T$ and the maps $f_i:=\bigcup_{0\leq j \leq i} (f_S)_j \times (f_T)_{i-j}\colon  (S\times T)_i:=\bigcup_{0\leq j \leq i} S_j \times T_{i-j}\to S\times T$ for $0\leq i \leq k+l$.
\end{Def}

The meaning of ``smooth'' in Definition \ref{defstrat} is only a special case of the general definition of smooth maps between the objects we introduced, and which will be given now.

\begin{Def} \label{defc}\label{defsmooth}
A continuous map \emph{$h\colon S\to M$} from an $n$-dimensional p-stratifold to a c-manifold is called \emph{smooth} if for all $j\leq n$ the composition $h f_j\colon S_j\to M$ is a c-map. 

A map \emph{$g\colon M\to S$} from a c-manifold to a p-stratifold is
called \emph{smooth} if $h g$ is smooth for all smooth maps $h\colon S\to \R$. 

A map $h\colon S\to T$ from an $n$-dimensional p-stratifold to another p-stratifold is called \emph{smooth} if for all $j\leq n$ and for every component $(S_j)_k$ of $S_j$ there is a c-map $h_{(j,k)}\colon (S_j)_k\to T_{l(j,k)}$ such that the composition $h f_j\colon (S_j)_k\to T$ can be factored as $(f_T)_{l(j,k)} h_{(j,k)}$.
\end{Def}

The definition of smooth maps between p-stratifolds in particular
implies that a diffeomorphism of p-stratifolds is a homeomorphism of
the underlying topological spaces preserving the p-stratifold
structure. 

\begin{Rmk}
A peculiarity of the above definition is that for a p-stratifold $S$
and a c-manifold $M$ one obtains the same definition of maps $f\colon
S\to M$ regardless of whether one considers $M$ as a c-manifold or as
a p-stratifold, but it makes a difference for maps $g\colon M\to
S$. This is of no importance in the following, the general definition
of maps between p-stratifolds being only included for completeness; we
will only need diffeomorphisms.
\end{Rmk}
 
Basic concepts of differential topology extend to p-stratifolds, such
as the existence of smooth partitions of unity and Brown's
theorem~\cite[Propositions~2.3 and~2.6]{kreck}. 

Every algebraic variety with isolated singularities admits a canonical
structure as a p-stratifold~\cite{grinberga}. More generally, every abstract
pre-stratified space in the sense of Mather~\cite{mather} is a cornered
p-stratifold, a generalization of a p-stratifold which we have not
defined here~\cite[Theorem~3.4.2]{grinbergb}. Furthermore, for any abstract
pre-stratified space $V$, there is a p-stratifold $S$ such that $V$ is
isomorphic to $S$ as a stratified
space~\cite[Corollary~3.4.4]{grinbergb}. Various classes of singular
spaces, such as Whitney stratified spaces and algebraic varieties, admit a
structure of an abstract pre-stratified space. So they admit a
structure of a p-stratifold.

Since we will define a bordism relation of p-stratifolds, we need to
introduce a notion of p-stratifolds with boundary.

\begin{Def} 
An $n$-dimensional \emph{p-stratifold with boundary} is given by a
triple $(S, \partial T, \varphi_T)$ as $T=S\cup_{\varphi_T} (\partial
T \times [0,\varepsilon))$, where $S$ is an $n$-dimensional
p-stratifold, $\partial T$ is an $(n-1)$-dimensional p-stratifold and
$\varphi_T\colon \partial T \times (0,\varepsilon )\to S$ is a
diffeomorphism onto its image. We require that $\varphi_T(\partial T \times
(0,\varepsilon/2])$ is closed in $S$ and that for any open subset
$U\subseteq \partial T$ satisfying the conditions of
Definition~\ref{defa} the image $\phi_T(U\times (0,\varepsilon))$ is
open in $S$.

The \emph{$i$-stratum} of a p-stratifold with boundary is given by
$T_i:=S_i\cup_{\varphi_T} ((\partial T)_{i-1}\times[0,\varepsilon))$. 
\end{Def}

One sees that the approach taken in this definition is analogous to an
approach one can take when looking at c-manifolds: One can view a
c-manifold $M$ as the union of its interior $\oM$ with $\partial
M\times [0,\varepsilon)$, glued together by a collar map $\varphi$
that is a representative of the given equivalence class of
collars. In the above definition, the p-stratifold $S$ plays the role
of $\oM$. The last condition ensures, loosely speaking, that $\varphi$
really glues $\partial T$ to an end of $S$  and not somewhere into its
interior. 

Analogously to Definitions~\ref{defa},~\ref{defb} and~\ref{defc} one
defines induced p-stratifold structures on open subsets, products and
smooth maps for p-stratifolds with boundary. 

\begin{Ex} 
The cone $CS=S\times [0,1]/S\times \{1\}$ is given by the triple $(\oC S, S, id)$. 
\end{Ex}

In this example, the cone point cannot be part of the boundary if
$\dim S \geq 1$, because then $pt \times (0,\varepsilon)$ is not open
in $\oC S$. 

We now define the notion of bordism of p-stratifolds. 

\begin{Def}
Two pairs $(S,g)$ and $(S',g')$ of $n$-dimensional p-stratifolds $S$
respectively $S'$ and continuous maps $g\colon S\to X$ respectively
$g'\colon S'\to X$ to a topological space $X$ are called
\emph{bordant} if there exists an $(n+1)$-dimensional p-stratifold $T$
with boundary such that $\partial T=S\amalg S'$ and a continuous map
$h\colon T\to X$ such that $h|_{\partial T}=g\amalg g'$.  
\end{Def}

We observe that the cylinder $S\times [0,1]$ gives a bordism of a
p-stratifold $S$ to itself and that two p-stratifolds can be glued
canonically along common boundary components thanks to the existence
of germs of collars. Thus ``bordism'' is a well defined equivalence
relation. 

For a topological space $X$ the equivalence classes $[S,g]$ of
$n$-dimensional closed p-stratifolds and continuous maps $g\colon S\to
X$ form a set~\cite[Proposition~4.1]{kreck}, and a group structure is
given by disjoint union. The neutral element is $[\emptyset]$, and any
element is of order maximal $2$, thus we obtain a $\Zz$-vector space.

We do not obtain an interesting theory when we try to define bordism
of all p-stratifolds, since for any p-stratifold $S$ the cone $CS$ is
a p-stratifold that gives a zero bordism. Thus we have to place
certain restrictions on the p-stratifolds ensuring that not all
p-stratifolds can be ``coned off'' -- that not for all allowed
p-stratifolds $S$ the cone $CS$ is again allowed. These restrictions
have to be compatible with the bordism relation. 

We call a class of p-stratifolds a bordism class if it has properties
allowing an equivalence relation ``bordism'' of its p-stratifolds and
if the restrictions are sufficient to ensure that this bordism theory
is a homology theory. 

\begin{Def}\label{defbordclass}
A class $\Cl$ of p-stratifolds is called a \emph{bordism class} if it
satisfies the following properties: 

\begin{enumerate}
\item If $S$ is a closed p-stratifold in $\Cl$, then $S\times I$ is in $\Cl$.

\item If $S$ and $S'$ are compact p-stratifolds in $\Cl$, then
      $S\amalg S'$ is in $\Cl$. If $T'$ and $T''$ are compact
      p-stratifolds with boundary such that $\partial T'=S\cup S'$ and
      $\partial T''=S'\cup S''$, with $S\cap S'=S'\cap S''=\partial
      S=\partial S'=\partial S''$ an $(n-2)$-dimensional p-stratifold,
      then $T'\cup_{S'}T''\in \Cl$, and
      $\partial(T'\cup_{S'}T'')=S\cup_{\partial S'}S''\in \Cl$.  

\item Given a compact $n$-dimensional p-stratifold $S\in \Cl$, if
      $\rho\colon S\to \R$ is a smooth map and $t\in \R$ is a regular
      value of $\rho$, then $\rho^{-1}(t)$ is in $\Cl$, and
      $\rho^{-1}([t,\infty))$ and $\rho^{-1}((-\infty,t])$ are in
      $\Cl$. 
\end{enumerate} 
\end{Def}

\begin{Prop}\label{homologietheoriebedingungen}
If $\Cl$ is a bordism class, then bordism of compact p-stratifolds in
$\Cl$ is well-defined and defines a homology theory $\H_{\Cl,*}$. 
\end{Prop}

\begin{proof}
Property 1 ensures reflexivity and Property 2 ensures transitivity of
the bordism relation. Symmetry is clear. Thus one can define bordism
groups $G_n(X)$ of $n$-dimensional closed p-stratifolds $S$ in $\Cl$
and continuous maps $g\colon S\to X$. Induced maps are given by
composition: A smooth map $f\colon X\to Y$ induces $f_*\colon
G_n(X)\to G_n(Y), [S,g]\mapsto [S,f g]$ for all~$n$.  

In order for $G_n(X)$ to define a homology theory, the Mayer-Vietoris
sequence is needed. For an open covering $U_1\cup U_2$ of $X$, a
boundary operator $d\colon G_n(U_1\cup U_2)\to G_{n-1}(U_1\cap U_2)$
in the Mayer-Vietoris sequence is obtained as
follows~\cite[Proposition~5.1]{kreck}.  

For $[S,g]\in G_n(U_1\cup U_2)$ we take a smooth map $\rho\colon S\to
\R$ such that $\rho(g^{-1}((U_1\cup U_2)\setminus U_1))=1$ and
$\rho(g^{-1}((U_1\cup U_2)\setminus U_2))=-1$. Now we choose a regular
value $t\in (-1,1)$ of $\rho$ and set $W:= \rho^{-1}(\{t\})$. One
shows that this is a p-stratifold of dimension $n-1$ by inductively
applying the analogous result for smooth manifolds. It is again
compact without boundary, and Property 3 ensures that
$[W,g|_W]=:d([S,g])$ is an element of $G_{n-1}(U_1\cap U_2)$.  

One uses Properties 2 and 3 to check that this is well defined and
leads to a long exact Mayer-Vietoris
sequence~\cite[Theorem~5.2]{kreck}.  
\end{proof}

Any homology theory thus obtained satisfies the disjoint union axiom
because we are using compact p-stratifolds: Let $I$ be any countable
index set, and let $X:=\amalg_{i\in I} X_i$ be the disjoint union of
topological spaces $X_i$. A compact
p-stratifold only has finitely many components, so any map from
a compact p-stratifold to $X$ only hits
finitely many components of the disjoint union. Thus we have 
\[
\H_{\Cl,*}(\amalg_{i\in I} X_i)\cong \bigoplus_{i\in I} \H_{\Cl,*}(X_i).
\]

One way to obtain interesting bordism classes is to restrict attention
to p-stratifolds whose strata of a fixed codimension are empty or
nonempty. The most prominent example is the special case ``codimension
1-stratum empty''. For all $n\geq 0$, one requires an $n$-dimensional
p-stratifold $S$ to have an empty $(n-1)$-stratum
$S_{n-1}=\emptyset$. For CW-complexes, the bordism theory defined this
way is isomorphic to singular homology with
$\Zz$-coefficients~\cite[Theorem~20.1]{kreck}. 

We sketch the basic ideas, since they are similar to the ones we apply
in the next section. It is easily seen that this restriction fulfills
the properties listed in Definition~\ref{defbordclass} and thus
defines a homology theory.  

The decisive step is to compute the coefficients of the theory. The
important observation here is that any p-stratifold of dimension $n
\geq 1$ can be ``coned off''. $CS$ is a zero bordism, an
$(n+1)$-dimensional p-stratifold with boundary $S$ fulfilling
$(CS)_{n}=\emptyset$ since $(CS)_n=S_{n-1}\times [0,1]=\emptyset$ for
all $n\geq 1$. But for $n=0$ we are in a different situation: The cone
over a $0$-dimensional p-stratifold $S$ is a $1$-dimensional
p-stratifold that does have a $0$-stratum, the cone point, thus it is
not an allowed zero bordism. We obtain a non-trivial $0$-th homology
group that is isomorphic to $\N_0(pt)=H_0(pt)$. The comparison
theorem is used to complete the argument. 

We now introduce a different restriction on p-stratifolds having to do
with the Euler characteristic of certain subspaces. In
section~\ref{sec3}, we show that this
restriction fulfills the properties listed in Definition~\ref{defbordclass}. Then we analyze
the homology theory thus obtained by computing its coefficients and
general structure. 

Given the $i$-skeleton $S^{(i)}:=\bigcup_{j\leq i}f_j(S_j)$ of an
$n$-dimensional p-stratifold $S$ ($i\leq n$) and representatives
$\varphi_l\colon \partial S_l\times [0,\varepsilon)\to S_l$ of the
germs of the collars of the $S_l$ for all $i+1\leq l \leq n$, one can
define an open neighborhood $W^{(i,n)}$ of $S^{(i)}$ inductively by
setting 
\begin{eqnarray*}
W^{(i,i)} & := & S^{(i)} \\
W^{(i,l)} & := & W^{(i,l-1)} \cup  f_l
(\varphi_l(f_l^{-1}(f_l(\partial S_l)\cap
W^{(i,l-1)})\times[0,\varepsilon))) \text{ for $i<l\leq n$}. 
\end{eqnarray*}
One defines a retract of this neighborhood onto $S^{(i)}$ by the
composition $r_i:=
(r_i)_{i,i+1} \ldots (r_i)_{n-1,n}$, where $(r_i)_{l,l+1}\colon
W^{(i,l+1)}\to W^{(i,l)}$ is given by 
\begin{eqnarray*}
(r_i)_{l,l+1} & = & \begin{cases}
  id_{W^{(i,l)}} & \text{on $W^{(i,l)}$} \\
  f_{l+1}\pr_1\varphi^{-1}_{l+1} f^{-1}_{l+1} & \text{on
    $W^{(i,l+1)}\setminus W^{(i,l)}$}.
  \end{cases}
\end{eqnarray*}
We call the map $r_i$ the retract given by the collars, since it is
basically just repeated projection along the collars onto the
boundaries of the strata.  

In the case of p-stratifolds $T$ with boundary, $T=(S,\partial T,
\varphi_T)$, the open neighborhood $W^{(i,n)}$ of the $i$-skeleton
$T_i$ is defined as 
\[
W^{(i,n)}:=W_{S}^{(i,n)}\cup_{\varphi_T} W_{\partial
  T}^{(i-1,n-1)}\times [0,\varepsilon),
\] 
with the retract given by the collars as above. One can also define
open neighborhoods of the interiors of the strata, 
\[
E^{(i,n)}:=r_{i}^{-1}(T_{i}\setminus (T_{i}\cap T^{(i-1)}))\subseteq
W^{(i,n)}.
\] 

\begin{Def}\label{defreg}
An $n$-dimensional p-stratifold $S$ is called \emph{regular} if for
all $0\leq i \leq n$ and for all $x\in f_i(\oS_i)$ there is an open
neighborhood $U_x$ of $x$ in $f_i(\oS_i)$ diffeomorphic to the open
$i$-ball $B^i$ such that there is a diffeomorphism of p-stratifolds
$\psi\colon  r_i^{-1}(U_x)\xrightarrow{\sim} B^{i}\times F$, with $F$
a p-stratifold whose $0$-stratum $F_0$ is a single point. 

A regular p-stratifold $S$ is called an \emph{Euler stratifold} if for
all $x\in S$ these $F$ have the property that the complement of the
$0$-stratum has even Euler characteristic: $\chi(F\setminus F_0)\equiv
0 \mod 2$. 

An $n$-dimensional p-stratifold $T$ with boundary is called regular if
its interior $S$ is regular, and Euler if its interior $S$ is Euler. 
\end{Def}

These are the objects used in the sequel. The regularity condition
means that the $E^{(i,n)}$ are locally trivial fiber bundles over
$S_{i}\setminus (S_{i}\cap S^{(i-1)})$, and we have $\pr_2 \psi(x) =
F_0$. For a p-stratifold $T$ with boundary the diffeomorphism
condition on $\varphi_T$ implies that with $S$ also $\partial T$ is
regular respectively Euler. 

\begin{Ex}\label{excman}
Any c-manifold  regarded as a p-stratifold is an Euler stratifold
since for all points $x$ we obtain $r_i^{-1}(U_x)=U_x$, thus $F=F_0$,
and the above conditions are trivially satisfied. 
\end{Ex}

\begin{Lem}
A product of regular p-stratifolds is again regular, and a product of Euler stratifolds is an Euler stratifold.
\end{Lem}

\begin{proof}
Let $S$ and $T$ be regular p-stratifolds. Then every $(s,t)\in
f_{i+j}(\oS_i\times\oT_j)\subseteq S\times T$ has an open neighborhood
$(r_S)_i^{-1}(U_s)\times (r_T)_j^{-1}(U_t) \cong B^i \times F\times
B^j \times G \cong (B^{i+j})\times (F\times G).$ For the fibers
$F\times G$ there is the open covering $(F\setminus F_0)\times G\cup
F\times (G\setminus G_0)$ of $(F\times G)\setminus (F_0 \times G_0)$,
the intersection being $(F\setminus F_0) \times (G\setminus G_0)$. If
$S$ and $T$ are Euler stratifolds, this leads to $\chi((F\times
G)\setminus (F_0 \times G_0))\equiv \chi((F\setminus F_0)\times G) +
\chi( F\times (G\setminus G_0)) + \chi((F\setminus F_0) \times
(G\setminus G_0))\equiv 0 \mod 2.$ 
\end{proof}

\section{Euler homology}\label{sec3}

In this section, we define and analyze Euler homology. The homology
groups will be obtained as bordism groups of Euler stratifolds. 

\begin{Thm}
Bordism groups $Eh_*(X)$ can be defined by 
$$ Eh_n(X):=\{ [S,g]\; |\; S \text{ $n$-dim.~closed Euler stratifold, } g\colon S\to X \text{ cont.} \}$$ 
for all $n\geq 0$ and yield a homology theory. 
\end{Thm}

\begin{proof}
We have to prove the properties of Definition~\ref{defbordclass}. By
Example~\ref{excman}, c-manifolds are Euler stratifolds, and the Euler
property is preserved by the product, thus for any Euler stratifold
$S$ also $S\times [0,1]$ is an Euler stratifold. Gluing of Euler
stratifolds along common boundary components leads to an Euler
stratifold. The Euler property is a local condition, thus we only have
to check whether it is true at the points where the Euler stratifolds
are glued together. There it follows from the fact that it holds on
the boundary and from the existence of collars. 

For the points $x$ in the preimage of a regular value $t$ of a smooth
map $\rho\colon T\to \R$ one applies the result for smooth manifolds
to the neighborhood $U_x\cong B^i$ of $x$ in $T^{(i)}\setminus
T^{(i-1)}$. One obtains a submanifold $V_x\subseteq U_x$ diffeomorphic
to $B^{i-1}$. Since there are collars with which $\rho$ is compatible
(see definitions \ref{defcmap} and \ref{defsmooth}), all points in the
neighborhood $r_i^{-1}(V_x)$ defined by these collars are in
$\rho^{-1}(\{t\})$. This neighborhood has the same fibers as the
original neighborhood of $x$ in $T$ (we could use the same collars to
define it), thus the Euler condition is again satisfied. 
\end{proof}

The central question now is to determine the coefficients $Eh_*(pt)$ of this homology theory. We first observe that any Euler stratifold that has even Euler characteristic is zero bordant since the cone over it is again an Euler stratifold.

\begin{Lem}\label{lemconeoff}
Let $S$ be an $n$-dimensional Euler stratifold. If $\chi(S)\equiv 0\mod 2$, then $[S]=0$ in $Eh_n(pt)$. If $\chi(S)\equiv 1 \mod 2$, then $[S]=[pt_n]$ in $Eh_n(pt)$.
\end{Lem}

\begin{proof}
It is clear that $CS$ is again a regular p-stratifold and that the Euler property is true for all $x$ outside the $0$-stratum. For the cone point $x=f_0((CS)_0)$ we observe that an open neighborhood $r_0^{-1}(x)=B^0\times F$ of $x$ in $CS$ is given by $V=S\times (1-\varepsilon, 1]/S\times \{1\}$, and so $\chi(F\setminus F_0)=\chi(S\times (1-\varepsilon, 1))=\chi(S)$. Thus if $\chi(S)\equiv 0 \mod 2$ then $CS$ is an Euler stratifold with boundary $S$, and therefore $[S]=0\in Eh_n(pt)$. If $\chi(S)\equiv 1 \mod 2$ then $\chi(S\cup pt_n)\equiv 0 \mod 2$ implying $[S\cup pt_n]=0$ and thus $[S]=[pt_n]\in Eh_n(pt)$.  
\end{proof}

Determining the coefficients of $Eh_n(pt)$ now boils down to the
question of whether or not $[pt_n]$ defines a nontrivial element in $Eh_n(pt)$. The following theorem states that the Euler characteristic of the boundary of a compact Euler stratifold is even, implying that the point cannot be zero bordant and is therefore a nontrivial generator of $Eh_n(pt)$ for every $n$. 

\begin{Thm} \label{chivonrandnull}
$\chi(\partial T) \equiv 0 \mod 2 $ for all compact Euler stratifolds $T$ with boundary.
\end{Thm}

\begin{proof}
Let $T$ be a compact $n$-dimensional Euler stratifold with boundary given by the triple $(S,\partial T, \varphi_T)$. The idea of the proof is to calculate the Euler characteristic $\chi(T)$ in two different ways, obtaining the equations 
\begin{equation} \label{gl1}
\chi(T)\equiv \sum_{i=0}^n \chi(T_i) + \chi(\partial T) \mod2 
\end{equation} 
and
\begin{equation} \label{gl2}
\chi(T)\equiv \sum_{i=0}^n \chi(T_i)  \mod2 .
\end{equation}
For compact smooth manifolds the theorem is known to be true, and the first equation is shown inductively using this information. The second equation explicitly uses the Euler property.

If $T$ is compact, then $\partial T$ is compact as well as all $T_i$, since $T_i$ is homeomorphic to $T^{(i)}\cap T\setminus W^{(i-1,n)}$, and thus a closed subset of a compact space.

\emph{Proof of equation \ref{gl1}:}

We cover $T$ by the open subsets $T_n$ and $W^{(n-1,n)}$. Their intersection $T_n\cap W^{(n-1,n)}$ is homotopy equivalent to $\partial(S_n)$. We also have the covering of $\partial (T_n)$ by open neighborhoods of $\partial(S_n)$ and $(\partial T)_{n-1}\times \{0\}$, their intersection being homotopy equivalent to $\partial ((\partial T)_{n-1})$. Thus we know that $ \chi(T)\equiv \chi(T_n) + \chi(W^{(n-1,n)}) + \chi (\partial (S_n)) \mod 2$ and $0\equiv \chi(\partial (T_n))\equiv \chi(\partial(S_n)) + \chi((\partial T)_{n-1}\times \{0\})+\chi(\partial((\partial T)_{n-1}))\mod 2.$ Noting that $W^{(n-1,n)}$ is homotopy equivalent to $T^{(n-1)}$ and that $\chi(\partial(T_n))\equiv 0\equiv \chi(\partial(\partial T)_{n-1}) \mod 2$, we obtain $\chi(T)\equiv\chi(T_n)+\chi(T^{(n-1)})+\chi((\partial T)_{n-1})\mod 2$. Since the equation is obviously true for $0$-dimensional Euler stratifolds, we inductively obtain $ \chi(T)\equiv \sum_{i=0}^n \chi(T_i) +  \sum_{i=1}^n \chi((\partial T)_{i-1}).$ Applying this equation to the $(n-1)$-dimensional compact p-stratifold (without boundary) $\partial T$, one obtains $\chi(\partial T)\equiv \sum_{i=0}^{n-1}\chi((\partial T)_i)$ and thus $\chi(T)\equiv \sum_{i=0}^n \chi(T_i) + \chi(\partial T).$

\emph{Proof of equation \ref{gl2}:}

We use the same procedure, but for a different open covering of $T$: This time $T=T\setminus T^{(0)}\cup E^{(0,n)}$, and we know that $r_0\colon  E^{(0,n)}\to T^{(0)}$ is a locally trivial fiber bundle. Thus the intersection $T\setminus T^{(0)}\cap E^{(0,n)}$ also is a locally trivial fiber bundle over $T^{(0)}$. The fibers $F$ of $r_0$ become $F\setminus F_0$, and of these we know that they have even Euler characteristic. Thus the total space of this bundle also has even Euler characteristic. We obtain $\chi(T)\equiv \chi(T\setminus T^{(0)}) + \chi( E^{(0,n)}) \mod 2$. $E^{(0,n)}$ is homotopy equivalent to $T_0$. To proceed inductively, we show that $\chi(T\setminus T^{(j)})=\chi(T\setminus T^{(j+1)}) + \chi(T_{j+1}).$ For this we use the covering $T\setminus T^{(j)}=T\setminus T^{(j+1)}\cup E^{(j+1,n)}$. By definition, the intersection of $E^{(j+1,n)}$ with $T\setminus T^{(j+1)}$ is a locally trivial fiber bundle and the fibers all have even Euler characteristic, so $\chi(T\setminus T^{(j)})=\chi(T\setminus T^{(j+1)})+ \chi( E^{(j+1,n)}) \mod 2.$ Since $ E^{(j+1,n)}\simeq \oT_{j+1}\simeq T_{j+1}$ and $T\setminus T^{(n)}=\emptyset$ we can proceed inductively to finally obtain equation \ref{gl2}. 
\end{proof}

Thus we have calculated the coefficients of $Eh_*$.
\begin{Cor}
$Eh_n(pt)\cong \Zz$ for all $n\geq 0$.
\end{Cor}

A natural transformation $\Phi\colon \N_*(X)\to Eh_*(X)$ is obtained by regarding an $n$-dimensional singular manifold $[M,f]\in\N_n(X)$ as an $n$-dimensional Euler stratifold $[M,f]\in Eh_n(X)$. Here $M$ is equipped with the canonical equivalence class of collars $\varphi\colon \emptyset \to M$. For $X=pt$ this map leads to 
\begin{eqnarray*}
\N_n(pt) & \to & Eh_n(pt) \cong \Zz\\
{}[M] & \mapsto & \chi(M) \mod 2.
\end{eqnarray*}
$\Phi$ is a well defined natural transformation since both theories are defined as bordism theories. A bordism in $\N_*(X)$ leads to a bordism in $Eh_*(X)$ by choice of an equivalence class of collars, and constructions like induced maps and the Mayer-Vietoris sequence are defined analogously.

The homology theory $Eh_*(X)$ has a product structure given by the product of Euler stratifolds, $Eh_n(X)\times Eh_l(Y)\to Eh_{n+l}(X\times Y)$. The above natural transformation thus induces an $\N_*$-module structure on $Eh_*(X)$. 

The next goal is to understand not only the coefficient groups $Eh_*(pt)$ of the Euler homology theory, but also the groups $Eh_*(X)$ for any CW-complex $X$. This is done by comparing it to homology theories we know, more precisely: An isomorphism 
$$\psi\colon  Eh_*(X)\xrightarrow{\sim} H_*(X;\Zz)\otimes_{\Zz} \Zz [t]$$ 
of homology theories is given.  

We equip $\Zz [t]$ with the $\N_*$-module structure induced by the map $\N_*\to \Zz[t], [N]\mapsto \chi(N)\cdot t^{\dim N}$. Then we can extend the natural transformation $\Phi$ to a natural transformation $\varphi\colon  \N_*(X)\otimes_{\N_*}\Zz[t]\to Eh_*(X)$ by defining it on $\Zz[t]$ as $\varphi(t^n):=[pt_n]$. By the calculations above this is compatible with the $\N_*$-module structure, thus well defined. The map thus obtained is an isomorphism:

\begin{Thm}
The map $\varphi\colon  \N_*(X)\otimes_{\N_*}\Zz[t]\to Eh_*(X), [M^i,f^i]\otimes t^{n-i} \mapsto [M^i_n,f^i]$ is a natural isomorphism of homology theories preserving the product and $\N_*$-module structures.
\end{Thm}

\begin{proof}
We first show that $\N_*(X)\otimes_{\N_*}\Zz[t]$ is a homology theory. Then we can use the comparison theorem to show that it is isomorphic to $Eh_*(X)$. 

We start with the natural isomorphism $\N_*(X)\cong H_*(X,\Zz)\otimes_{\Zz} \N_*$ of $\N_*$-modules following from the work of Thom (see~\cite{stong} or~\cite{quillen}). One can tensor both sides over $\N_*$ with $\Zz[t]$ to obtain $\N_*(X)\otimes_{\N_*}\Zz[t]\cong H_*(X,\Zz)\otimes_{\Zz} \Zz[t].$ The right side is a homology theory since $\Zz[t]$ is free over $\Zz$, so the left side is a homology theory, too.

The natural transformation $\varphi$ preserves the homology theory structure since induced maps and the boundary operator in the Mayer-Vietoris sequence are defined analogously. It is an isomorphism for $X=pt$ since in this case it reduces to $\Zz[t]\to Eh_*(pt)\cong \Zz[t], t^n\mapsto [pt_n].$ By the comparison theorem, $\varphi$ is a natural isomorphism of homology theories for CW-complexes.

The map preserves the product and $\N_*$-module structures since they are defined analogously on both sides.
\end{proof}

In the above proof we see that $\N_*(X)\otimes_{\N_*}\Zz[t]$, and thus also $Eh_*(X)$, is isomorphic to $H_*(X;\Zz) \otimes_{\Zz} \Zz[t]$. (An $\N_*$-module structure is induced on $H_*(X;\Zz) \otimes_{\Zz} \Zz[t]$ by the one we defined on $\Zz[t]$.) We can even give an isomorphism in a concrete way:

\begin{Thm}
There is an isomorphism of homology theories given by
\begin{eqnarray*}
\mu \colon   \N_*(X)\otimes_{\N_*}\Zz[t] & \xrightarrow{\sim} & H_*(X;\Zz) \otimes_{\Zz} \Zz[t]\\
  {}[M^i,f^i]\otimes t^{n-i} & \mapsto & \sum_{0 \leq j \leq i} f^i_* \PD w_{i-j}(M^i) \otimes t^{n-j}.
\end{eqnarray*}
Here $w_{i-j}(M^i)\in H^{i-j}(M^i;\Zz)$ are the Stiefel-Whitney classes. The map $\PD\colon  H^{i-j}(M^i;\Zz)\to  H_j(M^i;\Zz)$ is the Poincaré duality map. This isomorphism preserves the product and $\N_*$-module structures.
\end{Thm}

\begin{proof}
We first need to show that $\mu$ is well defined. For this we need to see that the map 
\[
\mu'\colon  \N_*(X)\to H_*(X;\Zz) \otimes_{\Zz} \Zz[t], [M^i,f^i] \mapsto \sum_{0 \leq j \leq i} f^i_* \PD w_{i-j}(M^i) \otimes t^{i-j}
\] 
is a well defined $\N_*$-linear map preserving the product structure. 

It is well defined if it is a bordism invariant, that is if $\mu'([M^i,f^i])=0$ for $(M^i,f^i)=\partial(W^{i+1},h^{i+1})$. If $(M^i,f^i)=\partial(W^{i+1},h^{i+1})$, the inclusion $i\colon M^i\to W^{i+1}$ induces an inclusion of the tangent bundles, and since $M^i=\partial W^{i+1}$ has a collar in $W^{i+1}$ naturality of the Stiefel-Whitney classes implies $w_{i-j}(M^i)=i^*(w_{i-j}(W^{i+1})) \mbox{ for all } j.$ Since $i_* \PD = \PD \delta^*,$ with $ \delta^*$ the boundary operator in cohomology of the exact sequence for the pair $(W^{i+1},M^i)$, we have $\mu'([M^i,f^i])= 0$: We know that $\delta^* i^*=0$, so for all $j$ 
\begin{eqnarray*}
f^i_* \PD w_{i-j}(M^i) & = & h^{i+1}_* i_* \PD w_{i-j}(M^i) = h^{i+1}_* \PD  \delta^* w_{i-j}(M^i) \\
& = & h^{i+1}_* \PD \delta^* i^* w_{i-j}(W^{i+1}) = 0.
\end{eqnarray*}

The map $\mu'$ is obviously additive. The Stiefel-Whitney classes and the Poincaré duality map are compatible with the product structure, so $\mu'$ is, too. 
%Habe ich nachgerechnet
It remains to show that the $\N_*$-module structure is preserved. We calculate
\begin{eqnarray*}
\lefteqn{\mu\bigl([N].[M^i,f^i]\bigr)}\\
 & = & \mu\bigl([N\times M, cst\times f^i]\bigr)\\
& = & \Bigl(\sum_{0 \leq l \leq \dim N} cst_* \PD w_{\dim N-l}(N)\times \sum_{0 \leq j \leq i} f^i_* \PD w_{i-j}(M^i)\Bigr)\otimes t^{\dim N-l+i-j}\\
& = & \Bigl(cst_* \PD w_{\dim N}(N)\times \sum_{0 \leq j \leq i} f^i_* \PD w_{i-j}(M^i)\Bigr)\otimes t^{\dim N}\cdot t^{i-j}\\
& = & \chi(N)\cdot\sum_{0 \leq j \leq i} f^i_* \PD w_{i-j}(M^i)\otimes t^{\dim N}\cdot t^{i-j} = [N].\mu\bigl([M^i,f^i]\bigr).
\end{eqnarray*}

Since $\mu$ is given by the identity on the second factor, these calculations suffice to show that $\mu$ is a well defined $\N_*$-bilinear map on the product, thus a well defined map on the tensor product.

Now we show that $\mu$ is a natural transformation of homology theories. Obviously $\mu$ commutes with induced maps. $\mu$ also commutes with the boundary operator $d$ in the Mayer-Vietoris sequence: If $[T,g]=d([M^i,f^i])$, then $T\subseteq M^i$ and $g=f^i|_T$ by construction. Let $i\colon T\to M^i$ denote the inclusion. The manifold $T$ has a bicollar in $M^i$, so naturality of the Stiefel-Whitney classes implies $w_{i-j}(T)=i^*(w_{i-j}(M^i))$ for all $j$. We know that $PD i^* = d PD \colon  H^{p}(T)\to H_{i-p-1}(M^i)$ for all $p$ by inspection of definitions. We calculate
\begin{eqnarray*}
\lefteqn{d\mu\bigl([M^i_n,f^i]\bigr)}\\
& = & \sum_{0 \leq j \leq i} d(f^i_* \PD w_{i-j}(M^i) \otimes t^{n-j})
=\sum_{0 \leq j \leq i} (f^i|_T)_* d \PD w_{i-j}(M^i) \otimes t^{n-j}\\
& = & \sum_{0 \leq j \leq i} (f^i|_T)_* \PD i^* w_{i-j}(M^i) \otimes t^{n-j}
 = \sum_{0 \leq j \leq i} (f^i|_T)_* \PD w_{i-j}(T) \otimes t^{n-j}\\
& = & \mu\bigl([T_{n-1},f^i|_T]\bigr)=\mu d\bigl([M^i_n,f^i]\bigr).
\end{eqnarray*}
Finally we show that $\mu$ is an isomorphism for $X=pt$. In this case
\begin{eqnarray*}
\mu\colon  \qquad Eh_*(pt) \to & H_*(pt;\Zz) \otimes_{\Zz} \Zz[t] & \cong \Zz[t]\\
{}[pt_n] \mapsto & cst_* \PD w_{0}(pt) \otimes t^{n} & = t^n \qquad \text{for all } n. \qedhere
\end{eqnarray*}
\end{proof}

The main result now is an easy corollary.

\begin{Thm}
There is a natural isomorphism of homology theories $$ Eh_*(X)\cong H_*(X;\Zz)\otimes_{\Zz} \Zz [t]$$ that preserves the product and $\N_*$-module structures. 
\end{Thm}
\begin{proof}
By composition with $\varphi^{-1}$ the map $\mu$ yields an isomorphism 
$$\psi:=\mu \varphi^{-1}\colon  Eh_*(X)\xrightarrow{\sim} H_*(X;\Zz)\otimes_{\Zz} \Zz [t].$$ 
Since  $\varphi$ and $\mu$ preserve the product and $\N_*$-module structures, $\psi$ does, too.
\end{proof}

This isomorphism has the property that $\psi \Phi\colon  \N_*(X)\to H_*(X;\Zz)\otimes_{\Zz} \Zz [t]$ maps $[M,f]\in \N_n(X)$ to $ \sum_{0 \leq j \leq n} f_* \PD w_{n-j}(M)\otimes t^{n-j}$ by construction, since $\psi \Phi=\mu \varphi^{-1}\Phi=\mu\varphi^{-1}\varphi|_{\N_*(X)}=\mu|_{\N_*(X)}$. 

We summarize the different natural transformations in the following diagram:
$$
\xymatrix{
\N_*(X) \ar[rd]^{\Phi} \ar[dd]_{id\otimes 1} & \\
& Eh_*(X) \ar[dd]^{\mu\varphi^{-1}}_{\cong}\\
\N_*(X)\otimes_{\N_*}\Zz[t] \ar[ur]^{\varphi}_{\cong} \ar[dr]^{\mu}_{\cong} & \\
 & H_*(X;\Zz) \otimes_{\Zz} \Zz[t] .
}$$

\section{Equivariant Euler homology}

Let $G$ be a discrete group. It is possible to define equivariant homology theories using bordism of p-stratifolds with a $G$-action. The discrete group $G$ is a $0$-dimensional smooth manifold, so we can endow $G\times S$ with the product p-stratifold structure. 

\begin{Def}
A p-stratifold $S$ is called a \emph{$G$-stratifold} if there is a smooth proper $G$-action on $S$, i.e., a smooth map of p-stratifolds 
$
\theta \colon  G\times S  \to  S,
 (g,s)  \mapsto  gs
$
such that
\begin{enumerate}
\item $es=s$
\item $(gh)s=g(hs)$
\item $G\times S  \to  S\times S,  (g,s)  \mapsto  (gs,s)$ is proper.
\end{enumerate}
\end{Def}

The map $\theta_g\colon S\to S, s\mapsto gs$ is in $Aut(S)$ for all $g\in G$. In particular, $g\oS_i=\oS_i$ for all $g\in G$, collars are preserved, collars of the strata as well as collars of $S$, and the collar maps are equivariant (with respect to the trivial $G$-action on the interval). Hence we can also define the operation of $G$ on each stratum separately (and not only on its interior $\oS_i \hookrightarrow S$) by continuing it onto the boundary.
Thereby, we could also have defined a $G$-stratifold by saying: The strata $S_i$ are $G$-manifolds with $G$-equivariant collars which are glued together with $G$-equivariant maps.

\begin{Def}
Let $X$ be a $G$-space. An $n$-dimensional \emph{singular $G$-stratifold} is a pair $(S,f)$, where $S$ is an $n$-dimensional cocompact $G$-stratifold and $f\colon S\to X$ is a $G$-equivariant map. 
\end{Def}

As for usual homology theories defined via p-stratifolds, we can
restrict ourselves to $G$-stratifolds satisfying certain conditions;
thereby we obtain many different $G$-homology theories depending on the conditions imposed. The construction of equivariant homology theories using $G$-stratifolds is analogous to the non-equivariant theory; technical details can be found in~\cite{weberequivstratis}. 

We call a class of $G$-stratifolds a bordism class if it fulfills the
equivariant analog of the properties stated in
Definition~\ref{defbordclass}. (The interval is endowed with the
trivial $G$-action, the boundary components in Property~2 have to be
$G$-invariant and $\rho\colon S\to \R$ in Property~3 has to be
$G$-invariant.) Such a class allows an equivalence relation
``bordism'' of its $G$-stratifolds and this bordism theory is a
homology theory.

These axioms are usually fulfilled for conditions given locally, i.e., in terms of a neighborhood of $x$ for all $x\in S$. If a regular p-stratifold $S$ is endowed with an action of $G$ making it into a $G$-stratifold, then for every $x$ one can choose the $U_x$ in Definition~\ref{defreg} such that the isotropy group $G_x$ of $G$ at $x$ acts on $U_x$. Since $G_x$ fixes $x$, it also operates on $F_x$, making $F_x$ into a $G_x$-stratifold whose $0$-stratum $F_{x,0}$ is a single point. (One can choose $U_x$ small enough such that $gU_x\cap U_x\neq \emptyset$ implies $g\in G_x$ since the discrete group $G$ acts properly. Then one can take the $G$-invariant open neighborhood $r_i^{-1}(GU_x)=Gr_i^{-1}(U_x)\cong G(B^i\times F_x)$ and take the connected component of $Gr_i^{-1}(U_x)$ containing $x$.) If we impose conditions on the fibers $F_x$, this gives a bordism class of $G$-stratifolds. 

If we restrict ourselves to a bordism class $\Cl$ of $G$-stratifolds with free proper $G$-operations, we obtain an explicit description of the usual $\Cl$-stra\-ti\-fold homology theory of the Borel construction $EG\times_G X$. We have an isomorphism $\H^{G, free}_{\Cl,*}(X)\cong \H_{\Cl,*}(EG\times_G X),$ where $\H_\Cl$ stands for the homology theory in question. Standard examples are $Eh$, Euler homology, and $H(-;\Zz)$ or $H(-;\Z)$, singular homology with $\Zz$- or $\Z$-coefficients. The construction is analogous to the one for bordism of smooth manifolds with free $G$-action, described for example in~\cite[page~50ff]{conner-floyd}. Details can be found in~\cite[Section~2]{weberequivstratis}.

For spaces $X$ that are proper $G$-CW-complexes, we even obtain an equivariant homology theory~\cite[Proposition~3.1]{weberequivstratis}. This means that we have an induction structure linking the homology theories associated to various $G$~\cite[page~198f]{lueck}. In particular, for an inclusion $\alpha\colon  H\hookrightarrow G$, the induction map yields an isomorphism
$$ \ind_\alpha\colon  \H^H_{\Cl,*}(pt)\xrightarrow{\sim} \H^G_{\Cl,*}(G/H).$$
This reduces the analysis of $\H^G_{\Cl,*}$ to the inspection of
$\H^H_{\Cl,n}(pt)$ for all finite $H\leq G$. Given these groups, one
can use the Mayer-Vietoris sequence to compute $\H^G_{\Cl,*}(X)$ for all proper $G$-CW-complexes $X$.

We can define \emph{equivariant Euler homology}, a $G$-homology theory
$Eh^G_*$, using $G$-stratifolds that are Euler stratifolds in the
non-equivariant sense. 

\begin{Prop} 
Let $H$ be a finite group. Then for all $n\geq 0$, 
\[
Eh^H_{n}(pt)\cong \Zz.
\]
\end{Prop}

\begin{proof}
Let $[S]\in Eh^H_{n}(pt)$. If $\chi(S)\equiv 0 \mod 2$, by Lemma~\ref{lemconeoff} we can cone off $S$, obtaining $[S]=[\partial CS]=0$. If $\chi(S)\equiv 1 \mod 2$ we endow the one-point space $pt$ with the trivial $H$-operation, making it an $n$-dimensional $H$-stratifold $pt_n$. Then we cone off $S\cup pt_n$, which is possible since $\chi(S\cup pt_n)\equiv 0 \mod 2$. We obtain $[S\cup pt_n]=[\partial(C(S\cup pt_n))]=0$, and so $[S]=[pt_n]$. Since $\chi(\partial S)=0$ for all Euler stratifolds $S$, by Theorem~\ref{chivonrandnull} we know that $pt_n$ cannot be zero bordant. Thus $Eh^H_n(pt)=\{0,[pt_n]\}\cong \Zz$. 
\end{proof}

\begin{Cor} Let $G$ be a discrete group, and let $H$ be a finite subgroup of $G$. Then $Eh^G_{n}(G/H)\cong \Zz$ for all $n\geq 0$.
\end{Cor}

This suffices to calculate $Eh^G_*(X)$ for all proper $G$-CW-complexes $X$.

\section{Advanced equivariant Euler homology}

A more advanced definition of $G$-Euler stratifold will now be introduced using the equivariant Euler characteristic with values in the Burnside ring. It leads to a homology theory $\widetilde{Eh}$ which we call \emph{advanced equivariant Euler homology}. In this section, we prove the following proposition.

\begin{Prop}\label{resultat} 
Let $G$ be a discrete group, and let $H$ be a finite subgroup of $G$. Then $$\widetilde{Eh}^G_{n}(G/H)\cong  V_H \quad \text{ for all }n\geq 0, $$
where $V_H$ is a $\Zz$-vector space with base $\{ H/K \}$, where $K$ belongs to a complete set of conjugacy class representatives of the collection of subgroups of $H$ having odd index in their normalizer.
\end{Prop}

For any finite group $H$, the $\Zz$-vector space $V_H$ can be seen as a $\Zz$-analog of the Burnside ring. Namely, if $\Omega_*$ denotes oriented bordism, $\Omega^H_0(pt;\Z)$ is known to be the Burnside ring $A(H)$~\cite{stong1970}. On the other hand, for unoriented bordism $\N_*$, $\N^H_0(pt)$ is the $\Zz$-vector space $V_H$ described above~\cite[Proposition~13.1]{stong1970}. Thereby, we have a ring homomorphism $A(H)\to V_H$. Additively, $A(H)=\bigoplus_{(K)\in c(H)}\Z \cdot [H/K]$ and $V_H=\bigoplus_{(K)\in \K } \Zz \cdot [H/K]$, with $c(H)$ the set of conjugacy classes $(K)$ of subgroups of $H$ and $\K=\{ (K)\in c(H) \; | \; |NK/K| \text{ odd}\}$.

Let $X$ be a finite $H$-CW-complex. The \emph{equivariant Euler characteristic with values in the Burnside Ring} of $X$ is defined by 
\[
\chi^H(X):=\sum_{(K)\in c(H)}\chi(H\setminus X_{(K)})\cdot [H/K] \in A(H).
\]
Here $X_{(K)}:=\{x\in X \, | \, H_x\in (K)\}$ is the set of points $x\in X$ with isotropy group conjugate to $K$. The projection $A(H)\to V_H, \chi^H(X)\mapsto \bar{\chi}^H(X)$ corresponds to the projection $\Z\to \Zz, \chi(X)\mapsto\bar{\chi}(X)$ in the non-equivariant case. 

Instead of placing restrictions on the Euler characteristic $\chi$ of $F_x\backslash F_{x,0}$, we place restrictions on the equivariant Euler characteristic $\chi^{G_x}$~\cite[Definition~6.1]{lueck-rosenberg} of $F_x\backslash F_{x,0}$, thereby obtaining a finer invariant. Since $G$ acts properly on a $G$-stratifold $S$, all isotropy groups $G_x$ for $x\in S$ are finite.

\begin{Def}
An $n$-dimensional regular $G$-stratifold $S$ is called a \emph{$G$-Euler stratifold} if for all $x\in S$ the above $F_x$ have the property that the complement of the $0$-stratum maps to $0$ in~$V_{G_x}$: $\bar{\chi}^{G_x}(F_x\setminus F_{x,0})= 0 \in V_{G_x}$.
\end{Def} 

An important property of a $G$-Euler stratifold $S$ is that all the fixed point sets $S^H$ for $H\leq G$ are Euler stratifolds in the non-equivariant sense. First we will show that the fixed point sets are again regular p-stratifolds.

\begin{Prop}\label{shstratifold}
Let $S$ be a regular $G$-stratifold, and let $H\leq G$. Then $S^H$ is a regular p-stratifold. 
\end{Prop}
\begin{proof}
We want to show that $T:=S^H$ is a regular p-stratifold. For a
$G$-manifold $M$, we know that $M^H$ is a collection of
submanifolds. So $(S_i)^H=\coprod_{j\leq i}S_{i,j}$, where $S_{i,j}$
is the collection of components of $(S_i)^H$ of dimension $j$, a
$j$-dimensional manifold. $S_{i,j}$ is a c-manifold since the action
of $G$ is collar-preserving. We define $T_j=\coprod_{j\leq i}
S_{i,j}$, i.e., we collect the strata of the same dimension.  

The map $(f_T)_j\colon T_j\to T^{(j-1)}$ has to be defined as
$(f_T)_j=\sum_{j\leq i}(f_S)_i|_{S_{i,j}}$. In order for this
definition to make sense we have to ensure that $(f_S)_i|_{S_{i,j}}$
lands in $T^{(j-1)}$. If $x\in \IIm((f_S)_i|_{S_{i,j}})$, this means
that $x\in S^{(i-1)}$, that $x=(f_S)_i(y)$ with $y\in\partial(S_i)^H$, and
that a neighborhood of $y$ in $(S_i)^H$ has dimension $j$. Since $x\in
S^{(i-1)}$, $x$ is in the interior of some $S_k$ with $k\leq i-1$. The
fact that $G$ operates collar-preserving means that there is an
$\varepsilon>0$ such that $g\varphi(z,t)=\varphi(gz,t)$ for all $g\in
G$, for all $z\in \partial S_i$ in a neighborhood of $y$ and for all
$t\leq \varepsilon$. This implies that for $y\in \partial S_i^H$ we
have $H\varphi(y,t)=\varphi(Hy,t)=\varphi(y,t)$ for $t\leq
\varepsilon$.  
                                                                           
$S$ is a regular p-stratifold. So there is an open neighborhood $U_x$ of $x$ in $f_k(\oS_k)$ diffeomorphic to the open $k$-ball $B^k$ such that there is a diffeomorphism of p-stratifolds $\psi\colon  r_k^{-1}(U_x)\xrightarrow{\sim} B^{k}\times F_x$, with $F_x$ a p-stratifold whose $0$-stratum $F_{x,0}$ is a single point. Here $r_k$ is the retract given by the collars. Set $V=r_k^{-1}(U_x)$. So we have         
\[
V^H\cong (B^{k}\times F_x)^H = (B^{k})^H\times F_x^H,
\]
where the last equation results from the fact that $G$ operates collar-pre\-serv\-ing. The notation $F_x^H$ makes sense since $G_x$ operates on $F_x$ and $H$ is a subgroup of $G_x$. We can choose the neighborhood small enough such that $(B^{k})^H\cong B^l$ for some $l\leq k$, and such that $F_x^H$ still has a $0$-stratum consisting of a single point. 

Choose a $\delta\leq \varepsilon$ such that $\varphi(y,\delta)$ is in this neighborhood. Since a neighborhood of $y$ in $(S_i)^H$ has dimension $j$, we know that $\varphi(y,\delta/2)$ has a neighborhood of dimension $j$ in $V^H$. So its image under the above isomorphism also has a neighborhood of dimension $j$ in $(B^k)^H\times F_x^H$, we call it~$W$. 

We know that $\pr_2(\varphi(y,\delta/2))\not\in F_{x,0}^H$, so the dimension of $\pr_2(W)$ is at least $1$. On the other hand, $\pr_1(W)$ is an open subset of $(B^k)^H\cong B^l$, thus of dimension $l$. We conclude that $j=\dim (W)=\dim (\pr_1(W))+\dim (\pr_2(W))\geq l+1$, which shows that $l\leq j-1$, so $x\in S_{k,l}\subseteq T^{(j-1)}$. We have shown that $(f_S)_i|_{S_{i,j}}$ indeed lands in $T^{(j-1)}$. 

During the proof we also saw that a neighborhood of $x\in S_{k,l}$ is given by $V^H\cong (B^{k}\times F_x)^H = (B^{k})^H\times F_x^H\cong B^l\times F_x^H$. So $T=S^H$ is a regular p-stratifold. 
\end{proof}

Now we show that the p-stratifolds $S^H$ are Euler stratifolds if $S$ is a $G$-Euler stratifold.

\begin{Prop}\label{letzte}
Let $S$ be a $G$-Euler stratifold, and $H\leq G$ a subgroup of $G$. Then $S^H$ is an Euler stratifold.
\end{Prop}
\begin{proof}
We have shown in Proposition \ref{shstratifold} that $S^H$ is a regular p-stratifold. We now need to see that $S^H$ has the Euler property: Let $x\in S^H\subseteq S$. Let $F_x$ be the fiber at $x$ in $S$, $F_{x,S^H}$ the fiber at $x$ in $S^H$. Then $F_{x,S^H}=F_x^H, F_{x,S^H,0}=F_{x,0}^H=pt$, and $F_x^H\backslash F_{x,0}^H= (F_x\backslash F_{x,0})^H$. Thus 
\[
\chi(F_{x,S^H}\backslash F_{x,S^H,0})=\chi((F_x\backslash F_{x,0})^H).
\]
Since $S$ is a $G$-Euler stratifold, we know that $\chi^{G_x}(F_x\backslash F_{x,0})=\sum a_K \cdot [G_x/K]$ with $a_K\cdot |N_{G_x}K/K|$ even. Thus $\chi((F_x\backslash F_{x,0})^H)=\sum a_K \cdot |(G_x/K) ^ H|.$ 
We know that $N_{G_x}K/K$ operates freely on $G_x/K$ by 
\begin{eqnarray*}
N_{G_x}K/K \times G_x/K & \to & G_x/K\\
(nK, gK) & \mapsto & g n^{-1} K.
\end{eqnarray*}
So it acts freely on any fixed point set. Thus $(G_x/K)^H$ is a disjoint union of orbits each of which is isomorphic to $N_{G_x}K/K$. So $|(G_x/K)^H|$ is divisible by $|N_{G_x}K/K|$, and $a_K \cdot |(G_x/K)^H|$ is divisible  by the even number $a_K \cdot |N_{G_x}K/K|$. Thus $\chi((F_x\backslash F_{x,0})^H)=\sum a_K \cdot |(G_x/K) ^ H|$ is even.
\end{proof}

Now we are ready to generalize the procedure employed in the construction of non-equivariant Euler homology. We know that $G$-Euler stratifolds form a bordism class since they are defined locally, only by imposing conditions on $F_x$. So the assignment 
\begin{eqnarray*}
\protect\widetilde{Eh}^G_n\colon \GTop & \to & \Ab \\
X & \mapsto & \{\; [S,f]\;|\; S \text{ $n$-dim.~closed~$G$-Euler str.}, f\colon S\to X \text{ $G$-equ.}\}\\
{}(g\colon X\to Y) & \mapsto & (g_*\colon [S,f]\mapsto [S,gf])
\end{eqnarray*}
for all $n\in \Nat$ ($\widetilde{Eh}^G_n:=0$ for $n< 0$) together with the usual boundary operator $d$ is a $G$-homology theory~\cite[Proposition 1.9]{weberequivstratis}. 

We calculate $\widetilde{Eh}^H_{*}(pt)$ for any finite group $H$. 

\begin{Prop}
Let $H$ be a finite group. Then $\widetilde{Eh}^H_{n}(pt)\cong V_H$ for all $n\geq 0.$
\end{Prop}

\begin{proof}
We set $\phi\colon  \widetilde{Eh}_n^H(pt)  \to  V_H, [S] \mapsto \bar{\chi}^H(S)$. The first task is to see that this assignment is well defined. We need to show that $\bar{\chi}^H(\partial S)=0\in V_H$ for all compact $H$-Euler stratifolds $S$ with boundary. We consider the following commutative diagram with exact rows, where all sums are taken over $c(H)$, i.e., over conjugacy classes of subgroups of $H$.
\[
\xymatrix@C=2pc{
 & & A(H) \ar[r]^{\pr} \ar@{^{(}->}[d]^{i} & V_H \ar@{^{(}->}[d]^{\bar{i}} \ar[r] & 0 \\
0 \ar[r] & \oplus_{(K)} \Z \cdot x_K \ar[r]^{\cdot 2} \ar[d]_{\cong}^{\psi} & \oplus_{(K)} \Z \cdot x_K \ar[r]^{\pr} \ar[d]_{\cong}^{\psi} & \oplus_{(K)} \Zz \cdot x_K \ar[r] \ar[d]_{\cong}^{\bar{\psi}} & 0 \\
0 \ar[r]  & \oplus_{(L)} \Z \cdot [H/L]  \ar[r]^{\cdot 2}  & \oplus_{(L)} \Z  \cdot [H/L] \ar[r]^{\pr}  & \oplus_{(L)} \Zz  \cdot [H/L] \ar[r] & 0 
}
\]

Here $x_K:={|NK/K|}^{-1}{[H/K]}$. The inclusion $i$ is given by $[H/K]\mapsto [H/K]=|NK/K|\cdot x_K$, and $\bar{i}$ is the inclusion given by $[H/K]\mapsto x_K$ for $|NK/K|$ odd. If $|NK/K|$ is odd, we have $\pr i ([H/K])= x_K=\bar{i} \pr([H/K])$. If $|NK/K|$ is even, we have $\pr i ([H/K])=0=\bar{i}\pr([H/K])$. So $\pr i = \bar{i} \pr$. 

The map $\psi$ is given by $x_K \mapsto {|NK/K|}^{-1}{ch([H/K])}$. Here $ch$ is the character map defined by $ch([H/K])=\{|(H/K)^L|\}_{(L)\in c(H)}$. It is injective and has the property that $ch(\chi^H(X))=\{\chi((X)^L)\}_{(L)\in c(H)}$ for any $H$-CW complex $X$~\cite[5.5.1 and 5.5.8]{tomdieck1979}~\cite[4.2]{lueck-rosenberg}. We know that $\psi$ is an isomorphism~\cite[Proposition~5.8.3]{tomdieck1979}. The five-lemma implies that the induced map $\bar{\psi}$ is also an isomorphism. By definition, we have $\psi i = ch$. 

Let $\chi^H(\partial S)=\sum_{(K)} a_K \cdot [H/K]\in A(H)$. We need to show that this maps to $0$ in $V_H$. By definition $ch(\chi^H(\partial S))=\{\chi((\partial S)^L)\}_{(L)\in c(H)}$. Since $(\partial S)^L=\partial (S^L)$, the fact that the $S^L$ are Euler stratifolds implies that $\chi((\partial S)^L)=\chi(\partial (S^L))\equiv 0 \mod 2$ (Theorem~\ref{chivonrandnull}). Thus $\pr ch(\chi^H(\partial S))=0$. We calculate
\[
\bar{i} \bar{\chi}^H(\partial S)
= \bar{i}\pr(\chi^H(\partial S))
= \pr i (\chi^H(\partial S))
 = \bar{\psi}^{-1} \pr \psi i (\chi^H(\partial S))=0.
\]
By injectivity of $\bar{i}$, this implies $\bar{\chi}^H(\partial S)=0$.

The map $\phi$ is injective: If $S$ is an $H$-Euler stratifold with $\phi([S])=\bar{\chi}^H(S)=0,$ then $CS$ is an $H$-Euler stratifold. The only point where there is something to check is the cone point. If $x$ is the cone point, $F_x\backslash F_{x,0}\simeq S$ and $H_x=H$, thus $\bar{\chi}^{H_x}(F_x\backslash F_{x,0})=\bar{\chi}^H(S)=0\in V_H$. So $[S]=[\partial CS]=0$.

The map $\phi$ is surjective since the $0$-dimensional $H$-manifolds $H/K$ can be seen as $n$-dimensional $H$-stratifolds $H/K_n$. Thus the elements $[H/K_n]$ with $|NK/K|$ odd are in $\widetilde{Eh}^H_{n}(pt)$ and their images form a basis of $V_H$. 
\end{proof}

Because of the induction structure, this suffices to obtain the result stated in Proposition~\ref{resultat} for all discrete groups $G$. Using the Mayer-Vietoris sequence, we can then calculate $\widetilde{Eh}^G_*(X)$ for all proper $G$-CW-complexes~$X$.

\bibliographystyle{alpha}

\end{document}